\documentclass[psamsfonts,reqno]{amsart}

\usepackage{amssymb}
\usepackage{amscd}
\usepackage{eucal}
\usepackage{psfig}
\usepackage[dvips]{graphicx}
\usepackage[small]{caption}

\newcommand{\beq}{\begin{equation}}
\newcommand{\eeq}{\end{equation}}
\newcommand{\sds}{\strut\displaystyle}

\newcommand{\Real}{\mbox{Re\,}}
\newcommand{\Imag}{\mbox{Im\,}}

\newcommand{\R}{\mathbb{R}}

\newcommand{\C}{\mathbb{C}}
\newcommand{\N}{\mathbb{N}}
\newcommand{\I}{\mathbb{I}}

\newcommand{\sfrac}[2]{{\vphantom1\smash{\lower.5ex\hbox{\small$#1$}}\over
        \vphantom1\smash{\raise.4ex\hbox{\small$#2$}}}} 

\newtheorem{theorem}{Theorem}
\newtheorem*{theorem-nonumber}{Theorem}

\newtheorem{lemma}{Lemma}

\theoremstyle{remark}
\newtheorem*{remark}{Remark}

\theoremstyle{definition}

\begin{document}

\title[Cauchy Integral Equations]{On the Solution of a Class of Cauchy Integral Equations}

\author[E. De Micheli]{E. ~De Micheli}
\address[E. De Micheli]{IBF -- Consiglio Nazionale delle Ricerche \\ Via De Marini, 6 - 16149 Genova, Italy}
\email[E.~De Micheli]{demicheli@ge.cnr.it}

\author[G. A. Viano]{G. A. ~Viano}
\address[G. A. ~Viano]{Dipartimento di Fisica - Universit\`a di Genova,
Istituto Nazionale di Fisica Nucleare - sez. di Genova\\
Via Dodecaneso, 33 - 16146 Genova, Italy}
\email[G.A.~Viano]{viano@ge.infn.it}

\begin{abstract}
In a previous paper we have presented a new method for solving a
class of Cauchy integral equations. In this work we discuss in
detail how to manage this method numerically, when only a finite
and noisy data set is available: particular attention is focused
on the question of the numerical stability.
\end{abstract}

\maketitle

\section{Introduction}
\label{introduction_section}
In a previous paper \cite{DeMicheli} we have proved the following theorems.

\begin{theorem}
\label{the:1}
Let us consider the following series
\beq
\label{uno}
\frac{1}{2\pi} \sum_{n=0}^\infty a_n z^n \qquad (z=x+iy;\, x,y \in \R),
\eeq
and suppose that the set of numbers $\{f_n\}_0^\infty$, $f_n=n^p a_n$ $(p\geq 0,~n=0,1,2,\ldots)$
satisfies the following Hausdorff conditions
\beq
\label{due}
(n+1)^{(1+\epsilon)} \sum_{i=0}^n \binom{\sds n}{\sds i}^{2+\epsilon}
\left | \Delta^i f_{(n-i)} \right |^{2+\epsilon} < M \qquad (n=0,1,2,\ldots; \epsilon > 0),
\eeq
where $M$ is a positive constant, and $\Delta$ is the
difference operator defined by: $\Delta f_n=f_{n+1}-f_n$,
$\Delta^k f_n= \sum_{m=0}^k (-1)^m \binom{k}{m}
f_{n+k-m}$. Then:
\begin{itemize}
\item[(1)] series $(\ref{uno})$ converges uniformly to a function $f(z)$ analytic in
the unit disk $D_0=\{z \mid |z| < 1\}$;
\item[(2)] $f(z)$ admits a holomorphic extension to the ``cut-plane''
$\{z \in \C \setminus (1, +\infty)\}$;
\item[(3)] the jump function $F(x)=-i[f_+(x)-f_-(x)]$,
$(f_\pm(x)=\lim_{\substack{\eta\rightarrow 0 \\ \eta > 0}} f(x\pm i\eta))$,
belongs to $L^2(1,+\infty)$, and, moreover, if $p\geq 1$, it is a function of class $C^{(p-1)}$.
\end{itemize}
\end{theorem}
The $\epsilon$ in (\ref{due}) (missing in Ref. \cite{DeMicheli}) is needed to guarantee the continuity of
$\tilde{a}(\lambda)$ (Carlsonian interpolation of the $a_n$'s) at $\lambda=-1/2+i\nu$, $\nu\in\R$
(see next formula (\ref{duenove})).

Let us remind the reader that these results in the case $p=0$ (i.e. $f_n=a_n$) are
due to Stein and Wainger \cite{Stein}; however, these authors do not formulate the requirements
on the Taylor coefficients $a_n$ by the use of the Hausdorff condition (\ref{due}).

\begin{theorem}
\label{the:2}
If in series $(\ref{uno})$ the coefficients $a_n$ satisfy the
Hausdorff condition $(\ref{due})$ (i.e., $f_n=a_n$, $n=0,1,2,\ldots$), then the jump
function $F(x)$ can be represented by the following expansion, that converges in
the sense of $L^2$-norm:
\beq
\label{tre}
F(x)=\sum_{m=0}^\infty c_m \phi_m(x) \qquad (x\in [1,+\infty)),
\eeq
where the coefficients $c_m$ are given by:
\beq
\label{quattro}
c_m=\sqrt{2} \sum_{n=0}^\infty \frac{(-1)^n}{n!}\, a_n
P_m \left [ -i \left (n+\frac{1}{2}\right )\right ],
\eeq
$P_m$ being the Pollaczek polynomials. The functions $\phi_m(x)$ form a basis in
$L^2(0,+\infty)$, and are expressed by:
\beq
\label{cinque}
\phi_m(x) = B_m(x) \frac{e^{-1/x}}{x},
\eeq
where $B_m(x)=i^m \sqrt{2} L_m(2/x)$, $L_m$ being the Laguerre polynomials.
\end{theorem}

From these results it derives that we can formally solve the following integral
equation of Cauchy type:
\beq
\label{sei}
f(z)=\int_1^{+\infty} \frac{F(x)}{x-z}\,dx,
\eeq
if the infinite set of Taylor coefficients $a_n=f^{(n)}(0)/n!$ $(n=0,1,2,\ldots)$ is
known, and, in addition, if they are supposed to satisfy the Hausdorff condition (\ref{due})
with $f_n=a_n$. But, at this point, a very serious problem is met: How can
this procedure be managed numerically? In fact, it must be noted that, in practice, only a finite number
of Taylor coefficients can be known, and, moreover, they are usually affected by numerical errors. Furthermore,
if these coefficients are the results of experimental measurements the situation is even worse
because also the error proper of any measurement must be considered. Finally,
let us remark that the determination of the jump function $F(x)$ starting from
an approximate knowledge of the function $f(z)$ is a typical example of
improperly posed problem in the sense of Hadamard: the solution does not depend continuously on
the data. In order to make this point more clear notice that the problem of solving the
Cauchy integral equation (\ref{sei}) is strictly connected to the problem of the
analytic continuation up to the cut. We can conformally map the cut $z$-plane
onto the unit disk in the $\zeta$-plane (i.e. the domain $|\zeta| < 1$); in this map
the upper (lower) lip of the cut is mapped in the upper (lower) half of the unit
circle. Therefore the continuation up to the cut corresponds to the continuation
up to the unit circle $|\zeta|=1$. As is well known, the uniqueness of the analytic
continuation does not guarantee its continuity in the $L^2$ or in the uniform topology.
Furthermore, performing the analytic continuation up to the boundary
of the analyticity domain (i.e. up to the unit circle in the $\zeta$-plane geometry)
is a severely ill-posed problem \cite{John}.

In the standard methods of regularization, whose credit is generally due to Tikhonov
\cite{Groetsch,Tikhonov}, normally some appropriate global bound on the solution are imposed, by assuming
some prior knowledge on the solution itself. In this way a subspace
of the solution space is determined, and then a solution belonging to this subspace is looked for.
If this subspace is compact, then continuity follows from compactness. In the numerical analytic
continuation in the unit disk ($\zeta$-plane geometry), if one avoids
going up to the boundary, this subspace is obtained by imposing a global bound
at the boundary. This point can be easily understood by observing
that from a bound on the functions at $|\zeta|=1$ a bound on their
derivatives can be derived (in any point inside the unit circle), via the Cauchy integral
representation. Conversely, if the analytic continuation is required
up to the cut (i.e. up to the boundary), then a global bound only on the
functions is not sufficient, and a uniform bound on the first derivative becomes necessary.
Then, in this case, the regularization by the use of a-priori bounds is rather cumbersome.
Furthermore, it must be observed that in several problems of physical interest the
prior knowledge which allows for imposing bounds on the solution is missing, or can be
very poor. See on this point Refs. \cite{Magnoli,Viano}, where the difference between
synthesis and/or inverse problems is widely discussed: while in the first
class of problems the prior knowledge is intrinsic to the problem itself, being a
part of its formulation, for the second class (i.e., inverse problems) this is
not true. We are thus motivated to go beyond the Tikhonov and Tikhonov--like methods, and,
accordingly, to look for regularization procedures that do not make use of prior knowledge.

The method that we present in this paper is very far from the Tikhonov--like regularization procedures, and
it is actually based on the following observation: when in expansion (\ref{tre}) we use the actual data set,
which is usually composed by only a finite number of data affected by noise, then, although
the formal series (\ref{tre}) diverges, nevertheless the effect of the errors
(and of the fact that the number of data is finite) remains quite small in the beginning of the
expansion, and there will exist a point (i.e. a certain value of $m$) where the divergence sets in.
Thus, the basic idea is to stop the expansion at the point where it turns to diverge.
This rough and qualitative description can be put in rigorous form by proving that even
though series (\ref{tre}) diverges, nevertheless it converges (in the sense of the
$L^2$-norm) as the number of data tends to infinity and the noise vanishes. This is the main point
of the paper and it will be proved in Section \ref{solution_section}.

A very delicate point of the method consists in determining a procedure for the determination
of the truncation point in the approximation which will be proved to converge
to the unknown jump function $F(x)$. This procedure, that will be illustrated in
Section \ref{numerical_section}, which is devoted to the numerical analysis, is based on the
asymptotic behavior of the Pollaczek polynomial $P_m[-i(n+1/2)]$ for large
values of $m$ (at fixed $n$), whose proof will be given in the Appendix.

\section{Solution of the Cauchy integral equation when the data are
a finite number of noisy Taylor coefficients}
\label{solution_section}

Let us assume that only a finite number of noisy Taylor coefficients $a_n^{(\epsilon)}$
are known, and that $|a_n^{(\epsilon)}-a_n|\leq\epsilon$
$(n=0,1,2,\ldots,N; \epsilon > 0)$. Let $c_m^{(\epsilon,N)}$ be the
following finite sums
\beq
\label{dueuno}
c_m^{(\epsilon,N)}=\sqrt{2} \sum_{n=0}^N \frac{(-1)^n}{n!}\, a_n^{(\epsilon)}
P_m\left [ -i\left (n+\frac{1}{2}\right )\right ].
\eeq
Accordingly, we have $c_m^{(0,\infty)}=c_m$ (see formula (\ref{quattro})).
We can then prove the following lemma.

\begin{lemma}
\label{le:1}
The following statements hold true:
\begin{itemize}
\item[(i)]
\beq
\label{duedue}
\sum_{m=0}^\infty \left | c_m^{(0,\infty)} \right |^2 = \|F\|_{L^2[1,+\infty)}^2 = C;
\eeq
\item[(ii)]
\beq
\label{duetre}
\sum_{m=0}^\infty \left | c_m^{(\epsilon,N)} \right |^2 = +\infty;
\eeq
\item[(iii)]
\beq
\label{duequattro}
\lim_{\substack{N \rightarrow\infty \\ \epsilon\rightarrow 0}}
c_m^{(\epsilon,N)} = c_m^{(0,\infty)} = c_m \qquad (m=0,1,2,\ldots);
\eeq
\item[(iv)] if $k_0(\epsilon,N)$ is defined as
\beq
\label{duecinque}
k_0(\epsilon,N)= \max \{k\in\N : \sum_{m=0}^k \left | c_m^{(\epsilon,N)}\right |^2 \leq C\},
\eeq
then
\beq
\label{duesei}
\lim_{\substack{N\rightarrow\infty \\ \epsilon\rightarrow 0}}
k_0(\epsilon,N) = + \infty;
\eeq
\item[(v)] the sum
\beq
\label{duesette}
M_k^{(\epsilon,N)}=\sum_{m=0}^k \left | c_m^{(\epsilon,N)}\right |^2 \qquad (k\in\N)
\eeq
satisfies the following properties:
\begin{enumerate}
\item[(a)] it increases for increasing values of $k$;
\item[(b)] the following relationships hold true
\beq
\label{dueotto}
M_k^{(\epsilon,N)} \geq \left |c_k^{(\epsilon,N)}\right |^2 \;\substack{\sim \\ k\rightarrow +\infty}\;
\frac{1}{(N!)^2} (2k)^{2N} \qquad (N ~\mbox{fixed}).
\eeq
\end{enumerate}
\end{itemize}
\end{lemma}

\begin{proof}
(i) In Ref. \cite{DeMicheli} we proved that the coefficients $a_n$ in expansion (\ref{uno}),
which are supposed to satisfy condition (\ref{due}), admit a unique interpolating function
$\tilde{a}(\lambda)$ ($\lambda=\sigma+i\nu$) which belongs to the Hardy space
$H^2(\C_{-1/2})$ with norm
$\|\tilde{a}\|_2=\sup_{\sigma>-1/2}(\int_{-\infty}^{+\infty}|\tilde{a}(\sigma+i\nu)|^2d\nu)^{1/2}$
($\C_{-1/2}=\{\lambda\in\C\,|\, \Real\lambda > -1/2\}$).
Furthermore, we proved that $\tilde{a}(-1/2+i\nu)$, i.e., the restriction of $\tilde{a}(\sigma+i\nu)$
to the line $\Real\lambda\equiv\sigma=-1/2$, belongs to $L^2(-\infty,+\infty)$, and can be represented
by the following expansion which converges in the sense of the $L^2$--norm:
\beq
\label{duenove}
\tilde{a}\left (-\frac{1}{2}+i\nu\right )=\sum_{m=0}^\infty d_m \psi_m (\nu),
\eeq
where $d_m=\sqrt{2\pi}c_m$ (see formula (\ref{quattro})), and the
functions $\psi_m(\nu)$ are the Pollaczek functions (see \cite{DeMicheli})
\beq
\label{polla}
\psi_m(\nu)=\frac{1}{\sqrt{\pi}} \Gamma\left (\frac{1}{2}+i\nu\right ) P_m(\nu),
\eeq
(see Refs. \cite{Bateman2,Szego} for what concerns the Pollaczek polynomials). In view of the Parseval equality
we have:
\beq
\label{duedieci}
\sum_{m=0}^\infty |d_m|^2 = \int_{-\infty}^{+\infty}
\left |\tilde{a}\left (-\frac{1}{2}+i\nu\right )\right |^2\, d\nu.
\eeq
Next, we recall that $\tilde{a}(-1/2+i\nu)$ is the Fourier transform of
$F(v) \exp(v/2)$, where $F(v)$ represents the jump function associated
with the cut located at $\theta=iv$ ($v\in (0,+\infty)$) in the complex $\theta$-plane geometry
(see Fig. 1A in Ref. \cite{DeMicheli} which refers to the series $(1/2\pi)\sum_{n=0}^\infty a_n e^{-in\theta}$).
Therefore, in view of the Plancherel theorem and equality
(\ref{duedieci}) we get (see formula (46) of Ref. \cite{DeMicheli}):
\beq
\label{dueundici}
\int_{-\infty}^{+\infty} \left |\tilde{a}\left (-\frac{1}{2}+i\nu\right )\right |^2\, d\nu =
2\pi \int_{0}^{+\infty} \left | e^{v/2} F(v)\right |^2\,dv = \sum_{m=0}^\infty |d_m|^2.
\eeq
Next we pass from the $\theta$-plane geometry of Theorem 2 of Ref. \cite{DeMicheli} to the $z$-plane geometry of
Theorem 2' of Ref. \cite{DeMicheli} through the transformation $z=\exp(-i\theta)$. Then, to the cut at
$\theta=iv$ there corresponds the cut in the $z$-plane geometry located on the real
$x$-axis ($x \equiv \Real z$) from $x=1$ up to $x=+\infty$. Accordingly, the jump function
associated with the cut will be $F(\ln x)$. Adopting the convention used in \cite{DeMicheli}
we still denote this latter jump function by $F(x)$ in order to avoid a useless
proliferation of symbols. Therefore, by the substitution: $e^v=x$, from
equality (\ref{dueundici}) we obtain:
\beq
\label{duedodici}
\frac{1}{2\pi}\sum_{m=0}^\infty |d_m|^2=\int_{1}^{+\infty} |F(x)|^2\,dx = C.
\eeq
Since $c_m^{(0,\infty)}=d_m/\sqrt{2\pi}$, from (\ref{duedodici}) equality (\ref{duedue}) follows. \\
(ii) Let us rewrite the sums $c_m^{(\epsilon,N)}$ as follows:
\beq
\label{duetredici}
c_m^{(\epsilon,N)}=\sum_{n=0}^N b_n^{(\epsilon)} P_m\left [-i\left (n+\frac{1}{2}\right )\right ],
\eeq
where $b_n^{(\epsilon)}=\sqrt{2}(-1)^n a_n^{(\epsilon)}/n!$. Now, we can
write the following inequality:
\beq
\label{duequattordici}
\begin{split}
\left | c_m^{(\epsilon,N)} \right | &=
\left | \sum_{n=0}^N b_n^{(\epsilon)} P_m\left [-i\left (n+\frac{1}{2}\right )\right ]\right | \\
& \geq \left | b_N^{(\epsilon)} P_m\left [-i\left (N+\frac{1}{2}\right )\right ]\right |
\cdot \left | 1-\frac{\left | \sum_{n=0}^{N-1} b_n^{(\epsilon)} P_m\left [-i\left (n+\frac{1}{2}\right )\right ]\right |}
{\left | b_N^{(\epsilon)} P_m\left [-i\left (N+\frac{1}{2}\right )\right ]\right |} \right |.
\end{split}
\eeq
In the lemma proved in the Appendix we show that the asymptotic behavior of the Pollaczek polynomials
$P_m[-i(n+1/2)]$ for large values of $m$ (at fixed $n$) is given by:
\beq
\label{duequindici}
P_m\left [-i\left (n+\frac{1}{2}\right )\right ] \;\substack{\sim \\ m\rightarrow\infty}\;
\frac{(-1)^m i^m}{n!} (2m)^n.
\eeq
Therefore, we have:
\beq
\label{duesedici}
\begin{split}
&\frac{\left | \sum_{n=0}^{N-1} b_n^{(\epsilon)} P_m\left [-i\left (n+\frac{1}{2}\right )\right ]\right |}
{\left | b_N^{(\epsilon)} P_m\left [-i\left (N+\frac{1}{2}\right )\right ]\right |} \\
&\qquad\leq \frac{\sum_{n=0}^{N-1}\left | b_n^{(\epsilon)} P_m\left [-i\left (n+\frac{1}{2}\right )\right ]\right |}
{\left | b_N^{(\epsilon)} P_m\left [-i\left (N+\frac{1}{2}\right )\right ]\right |}
\;\substack{\sim \\ m\rightarrow\infty}\;
\sum_{n=0}^{N-1}\left | \frac{b_n^{(\epsilon)}}{b_N^{(\epsilon)}} \right | \frac{N!}{n!} (2m)^{n-N}
\xrightarrow[m\rightarrow\infty]{} 0.
\end{split}
\eeq
From (\ref{duequattordici}), (\ref{duequindici}) and (\ref{duesedici}) it follows that for $m$
sufficiently large:
\beq
\label{duediciassette}
\left |c_m^{(\epsilon,N)}\right |
\substack{\sim \\ m\rightarrow\infty} \frac{\left |b_N^{(\epsilon)}\right |}{N!} (2m)^N.
\eeq
Therefore, $\lim_{m\rightarrow\infty}|c_m^{(\epsilon,N)}|=+\infty$, and statement (ii) follows. \\
(iii) We can write the difference $[c_m^{(0,\infty)}-c_m^{(\epsilon,N)}]$ as
\beq
\label{duediciotto}
\begin{split}
c_m^{(0,\infty)}-c_m^{(\epsilon,N)}&= \sqrt{2}\left\{
\sum_{n=0}^N \frac{(-1)^n}{n!} (a_n-a_n^{(\epsilon)})
P_m\left [-i\left (n+\frac{1}{2}\right )\right ]\right . \\
&\left . + \sum_{n=N+1}^\infty \frac{(-1)^n}{n!} a_n P_m\left [-i\left (n+\frac{1}{2}\right )\right ]\right\}.
\end{split}
\eeq
In view of the fact that $\sqrt{2} \sum_{n=0}^\infty \frac{(-1)^n}{n!} a_n P_m[-i(n+\frac{1}{2})]$
converges to $c_m^{(0,\infty)}$, it follows that the second term in bracket (\ref{duediciotto})
tends to zero as $N\rightarrow \infty$. Concerning the first term we may write the inequality:
\beq
\label{duediciannove}
\left | \sum_{n=0}^N \frac{(-1)^n}{n!}(a_n-a_n^{(\epsilon)})P_m\left [-i\left (n+\frac{1}{2}\right )\right ]\right |
\leq \epsilon \sum_{n=0}^N \frac{1}{n!}\left | P_m\left [-i\left (n+\frac{1}{2}\right )\right ]\right |,
\eeq
where the inequalities $|a_n-a_n^{(\epsilon)}|\leq\epsilon$ $(n=0,1,2,\ldots,N)$ have been used.
Next, by rewriting the Pollaczek polynomials $P_m[-i(n+1/2)]$ as
\beq
\label{dueventi}
P_m\left [-i\left (n+\frac{1}{2}\right )\right ] = \sum_{j=0}^m p_j^{(m)} \left (n+\frac{1}{2}\right )^j,
\eeq
and substituting this expression in inequality (\ref{duediciannove}) we obtain
\beq
\label{dueventuno}
\epsilon \sum_{n=0}^N \frac{1}{n!}\left [\sum_{j=0}^m \left | p_j^{(m)}\right | \left (n+\frac{1}{2}\right )^j\right ].
\eeq
Next, we perform the limit for $N\rightarrow\infty$. In view of the fact that
$\sum_{j=0}^m p_j^{(m)} (n+1/2)^j$ is finite, and the series $\sum_{n=0}^\infty \frac{(n+1/2)^j}{n!}$
converges, we can exchange the order of the sums and write:
\beq
\label{dueventidue}
\epsilon \sum_{j=0}^m \left | p_j^{(m)}\right | \sum_{n=0}^\infty \frac{1}{n!} \left (n+\frac{1}{2}\right )^j.
\eeq
Finally, performing the limit for $\epsilon\rightarrow 0$, and recalling equality
(\ref{duediciotto}), statement (iii) is obtained. \\
(iv) From definition (\ref{duecinque}) it follows that for $k_1=k_0+1$, we have:
\beq
\label{dueventitre}
\sum_{m=0}^{k_1} \left | c_m^{(\epsilon,N)}\right |^2 > C.
\eeq
Statement (iv) (formula (\ref{duesei})) is proved if we can show that
$\lim_{\substack{N\rightarrow\infty \\ \epsilon\rightarrow 0}} k_1(\epsilon,N)=+\infty$.
Let us suppose that
$\lim_{\substack{N\rightarrow\infty \\ \epsilon\rightarrow 0}} k_1(\epsilon,N)$ is finite.
Then there should exist a finite number  $K$ (independent of $\epsilon$ and $N$), such that
$\lim_{\substack{N\rightarrow\infty \\ \epsilon\rightarrow 0}} k_1(\epsilon,N)\leq K$.
Then from inequality (\ref{dueventitre}) we have:
\beq
\label{dueventiquattro}
C < \sum_{m=0}^{k_1(\epsilon,N)}\left |c_m^{(\epsilon,N)}\right |^2\leq
\sum_{m=0}^{K} \left |c_m^{(\epsilon,N)}\right |^2.
\eeq
But as $N\rightarrow\infty, \epsilon\rightarrow 0$ we have (recalling also statement
(iii) formula (\ref{duequattro}))
\beq
\label{dueventicinque}
C < \sum_{m=0}^{K} \left |c_m^{(0,\infty)}\right |^2 < \sum_{m=0}^{\infty}\left |c_m^{(0,\infty)}\right |^2=C,
\eeq
which leads to a contradiction. Then statement (iv) follows.

(v) Concerning statement (a), it follows obviously from definition (\ref{duesette}) of
$M_k^{(\epsilon,N)}$.
Finally, the first one of relationships (\ref{dueotto}) is obvious; the second one follows
from the asymptotic behavior of $P_m[-i(n+1/2)]$ at large $m$ (for fixed $n$), proved in the
lemma of the Appendix.
\end{proof}

\begin{remark}
From statement (v) and formula (\ref{duesei}) it follows that,
for large values of $N$ and small values of $\epsilon$, the sum $M_k^{(\epsilon,N)}$ presents
a plateau in a neighborhood of $k=k_0$, if the function $F(x)$ being reconstructed is
sufficiently regular (see Section \ref{numerical_section} and figs. \ref{figura_1}B,
\ref{figura_2}B, \ref{figura_3}A, \ref{figura_4}A).
\end{remark}

Now we may introduce the following approximation of the jump function $F(x)$ (see expansion (\ref{tre})):
\beq
\label{dueventisei}
F^{(\epsilon,N)}(x)=\sum_{m=0}^{k_0(\epsilon,N)} c_m^{(\epsilon,N)}\psi_m(x) \qquad (x \in (1,+\infty)).
\eeq
Approximation $F^{(\epsilon,N)}(x)$ is defined through the truncation number
$k_0(\epsilon,N)$; the latter can be numerically determined by plotting the sum $M_k^{(\epsilon,N)}$
versus $k$, and exploiting properties (a) and (b), proved in statement
(v) of the previous lemma, and the property stated in the remark above.
For a more detailed analysis of this point see Section \ref{numerical_section}
where some numerical examples are given.

Now, we want to prove that the approximation $F^{(\epsilon,N)}(x)$ converges asymptotically to $F(x)$
in the sense of the $L^2$-norm, as $N\rightarrow\infty$ and $\epsilon\rightarrow 0$. We can prove
the following theorem.

\begin{theorem}
\label{the:3}
The equality
\beq
\label{dueventisette}
\lim_{\substack{N\rightarrow\infty \\ \epsilon\rightarrow 0}}
\left \| F - F^{(\epsilon,N)} \right \|_{L^2[1,+\infty)} = 0
\eeq
holds true.
\end{theorem}

\begin{proof}
From the Parseval equality it follows that:
\beq
\label{dueventotto}
\left \| F - F^{(\epsilon,N)} \right \|_{L^2[1,+\infty)}^2 =
\left\{\sum_{m=k_0+1}^{\infty}\left |c_m^{(0,\infty)}\right |^2 +
\sum_{m=0}^{k_0} \left |c_m^{(\epsilon,N)}-c_m^{(0,\infty)}\right |^2\right\}.
\eeq
Since $\sum_{m=0}^\infty\left |c_m^{(0,\infty)}\right |^2=C$ and
$\lim_{\substack{N\rightarrow\infty \\ \epsilon\rightarrow 0}} k_0(\epsilon,N)=+\infty$,
it follows that\\
$\lim_{\substack{N\rightarrow\infty \\ \epsilon\rightarrow 0}}
\sum_{m=k_0+1}^{\infty} |c_m^{(\epsilon,N)}|^2 = 0.$

It is then convenient to rewrite the second term of the r.h.s. of formula (\ref{dueventotto})
as follows. Let us denote by
\beq
\label{g0}
g_m^{(0,\infty)} =
\begin{cases}
c_m^{(0,\infty)} & \text{if $m$ is even}; \\
-i c_m^{(0,\infty)} & \text{if $m$ is odd};
\end{cases}
\eeq
\beq
\label{g1}
g_m^{(\epsilon,N)} =
\begin{cases}
c_m^{(\epsilon,N)} & \text{if $m$ is even}; \\
-i c_m^{(\epsilon,N)} & \text{if $m$ is odd}.
\end{cases}
\eeq
Notice that $\sum_{m=0}^{k_0}\left |c_m^{(\epsilon,N)}-c_m^{(0,\infty)}\right |^2 =
\sum_{m=0}^{k_0}\left (g_m^{(\epsilon,N)}-g_m^{(0,\infty)}\right )^2$,
and $g_m^{(0,\infty)}$ and $g_m^{(\epsilon,N)}$ are real. Next, we
introduce the following functions:
\begin{subequations}
\label{dueventinove}
\begin{eqnarray}
G^{(0,\infty)}(x)&=&\sum_{m=0}^{\infty}g_m^{(0,\infty)} \I_{[m,m+1[}(x),
\label{dueventinove.a} \\
G^{(\epsilon,N)}(x)&=&\sum_{m=0}^{\infty}g_m^{(\epsilon,N)} \I_{[m,m+1[}(x),
\label{dueventinove.b}
\end{eqnarray}
\end{subequations}
where $\I_E$ is the characteristic function of the set $E$. From statements (i),
(ii), and (iii) of the previous lemma (formulae (\ref{duedue}), (\ref{duetre}) and
(\ref{duequattro})) we obtain
\beq
\label{duetrenta}
\int_0^{+\infty} \left (G^{(0,\infty)}(x)\right )^2\, dx =
\sum_{m=0}^\infty \left (g_m^{(0,\infty)}\right )^2 = C,
\eeq
\beq
\label{duetrentuno}
\int_0^{+\infty} \left (G^{(\epsilon,N)}(x)\right )^2\, dx =
\sum_{m=0}^\infty \left (g_m^{(\epsilon,N)}\right )^2 = +\infty,
\eeq
\beq
\label{duetrentadue}
G^{(\epsilon,N)}(x) \xrightarrow[\substack{N\rightarrow\infty \\ \epsilon\rightarrow 0}]{}
G^{(0,\infty)}(x) \qquad (x\in [0,+\infty)).
\eeq
Hereafter we assume, for the sake of simplicity and without loss of
generality, that every term $g_m^{(\epsilon,N)}$ is different from zero. Next, let
$X(\epsilon,N)$ be the unique root of the equation
$\int_0^X\left (G^{(\epsilon,N)}(x)\right )^2\,dx=C$. Let us indeed observe that
$\int_0^X\left (G^{(\epsilon,N)}(x)\right )^2\,dx$ is a continuous non--decreasing function which is
zero for $X=0$, and $+\infty$ for $X\rightarrow +\infty$. Furthermore, from statement (iv) of the
previous lemma (formula (\ref{duesei})) we have
$\lim_{\substack{N\rightarrow\infty \\ \epsilon\rightarrow 0}} X(\epsilon,N)=+\infty$. \\
Then we can write
\beq
\label{duetrentatre}
\begin{split}
&\int_{0}^{X(\epsilon,N)}\left [G^{(\epsilon,N)}(x)-G^{(0,\infty)}(x)\right ]^2\,dx \\
&=\int_{X(\epsilon,N)}^{+\infty} \left (G^{(0,\infty)}(x)\right )^2dx
-2\int_{0}^{X(\epsilon,N)} G^{(0,\infty)}(x)
\left [ G^{(\epsilon,N)}(x)-G^{(0,\infty)}(x)\right ]\, dx.
\end{split}
\eeq
Next, we perform the limit for $N\rightarrow\infty$ and $\epsilon\rightarrow 0$.
Concerning the first term at the r.h.s. of formula (\ref{duetrentatre}) we have
\beq
\label{duetrentaquattro}
\lim_{\substack{N\rightarrow\infty \\ \epsilon\rightarrow 0}}
\int_{X(\epsilon,N)}^{+\infty} \left (G^{(0,\infty)}(x)\right )^2\,dx = 0.
\eeq
For what concerns the second term we introduce the following function:
\beq
\label{duetrentacinque}
H^{(\epsilon,N)}(x)=
\begin{cases}
G^{(\epsilon,N)}(x)-G^{(0,\infty)}(x) & \text{if $0\leq x\leq X(\epsilon,N)$}, \\
0 & \text{if $x > X(\epsilon,N)$}.
\end{cases}
\eeq
Then, we have by the use of the Schwarz inequality
\beq
\label{duetrentasei}
\int_{0}^{+\infty} \left | H^{(\epsilon,N)}(x)\right |^2\, dx \leq 4C \qquad (N < \infty, \epsilon > 0).
\eeq
Moreover, from (\ref{duetrentadue}) we have
\beq
\label{duetrentasette}
H^{(\epsilon,N)}(x)\xrightarrow[\substack{N\rightarrow\infty \\ \epsilon\rightarrow 0}]{}
0 \qquad x \in [0,+\infty).
\eeq
The family of functions $\{H^{(\epsilon,N)}\}$ is bounded in $L^2[0,+\infty)$, therefore it has
a subsequence which is weakly convergent in $L^2[0,+\infty)$.
The limit of this subsequence is zero. In fact, let us observe that
$|H^{(\epsilon,N)}(x)| \leq 2C$; then we consider the function
$H^{(\epsilon,N)}(x)\Phi(x)$, where $\Phi(x)$ is an arbitrary element of the class
of functions $C^\infty_c(\R^+)$. Then, we have $|H^{(\epsilon,N)}(x)\Phi(x)| \leq 2C|\Phi(x)|$,
and this inequality does not depend on $N$ and $\epsilon$.
In view of the Lebesgue dominated convergence theorem we can then write (see also the limit
(\ref{duetrentasette}))
\beq
\label{duetrentotto}
\lim_{\substack{N\rightarrow\infty \\ \epsilon\rightarrow 0}} \sup
\left | \int_0^{+\infty} H^{(\epsilon,N)}(x)\Phi(x)\,dx \right | = 0.
\eeq
Since the set of functions $C^\infty_c(\R^+)$ is everywhere dense in $L^2[0,+\infty)$,
given an arbitrary function $\psi \in L^2[0,+\infty)$ and an arbitrary number $\eta > 0$,
there exists a function $\Phi_k \in C^\infty_c(\R^+)$ such that
$\|\psi - \Phi_k \|_{L^2[0,+\infty)} < \eta$. Furthermore, through the Schwarz inequality
we have
\beq
\label{duetrentottobis}
\begin{split}
&\int_0^{+\infty}\left|H^{(\epsilon,N)}(x)\left[\psi(x)-\Phi_k(x)\right]\right|\,dx\\
&\qquad\leq \left(\int_0^{+\infty}\left| H^{(\epsilon,N)}(x)\right|^2\right)^{1/2}
\left(\int_0^{+\infty}\left|\psi(x)-\Phi_k(x)\right|^2\right)^{1/2}\leq 2\sqrt{C}\eta.
\end{split}
\eeq
From equality (\ref{duetrentotto}) and inequalities (\ref{duetrentottobis}) we thus conclude that
\beq
\label{duetrentottotris}
\lim_{\substack{N\rightarrow\infty \\ \epsilon\rightarrow 0}} \sup
\left | \int_0^{+\infty} H^{(\epsilon,N)}(x)\psi(x)\,dx \right | = 0,
\eeq
for any $\psi\in L^2[0,+\infty)$.

Next, by using the same type of arguments,  we can state that if there is an arbitrary subsequence
belonging to the family $\{H^{(\epsilon,N)}\}$, that converges weakly in $L^2[0,+\infty)$, then
the weak limit of this subsequence is necessarily zero. Finally, from the uniqueness of the
(weak) limit point it follows that the whole family $\{H^{(\epsilon,N)}\}$ converges weakly to
zero in $L^2[0,+\infty)$. We can thus write
\beq
\label{duetrentanove}
\lim_{\substack{N\rightarrow\infty \\ \epsilon\rightarrow 0}}
\int_0^{+\infty} G^{(0,\infty)}(x) H^{(\epsilon,N)}(x)\, dx = 0,
\eeq
and from inequality (\ref{duetrentatre}) it follows
\beq
\label{duequaranta}
\lim_{\substack{N\rightarrow\infty \\ \epsilon\rightarrow 0}}
\int_0^{X(\epsilon,N)} \left [G^{(\epsilon,N)}(x)-G^{(0,\infty)}(x)\right ]^2\, dx = 0.
\eeq
Since $\sum_{m=0}^{k_0}\left |c_m^{(\epsilon,N)}-c_m^{(0,\infty)}\right |^2
\leq \int_0^{X(\epsilon,N)} \left [G^{(\epsilon,N)}(x)-G^{(0,\infty)}(x)\right ]^2\, dx$,
we have
\beq
\label{duequarantuno}
\lim_{\substack{N\rightarrow\infty \\ \epsilon\rightarrow 0}}
\sum_{m=0}^{k_0}\left |c_m^{(\epsilon,N)}-c_m^{(0,\infty)}\right |^2 = 0,
\eeq
and, in view of equality (\ref{dueventotto}), the theorem is proved.
\end{proof}

\section{Numerical analysis}
\label{numerical_section}
\begin{figure}[ht]
\begin{center}
\includegraphics[width=11cm]{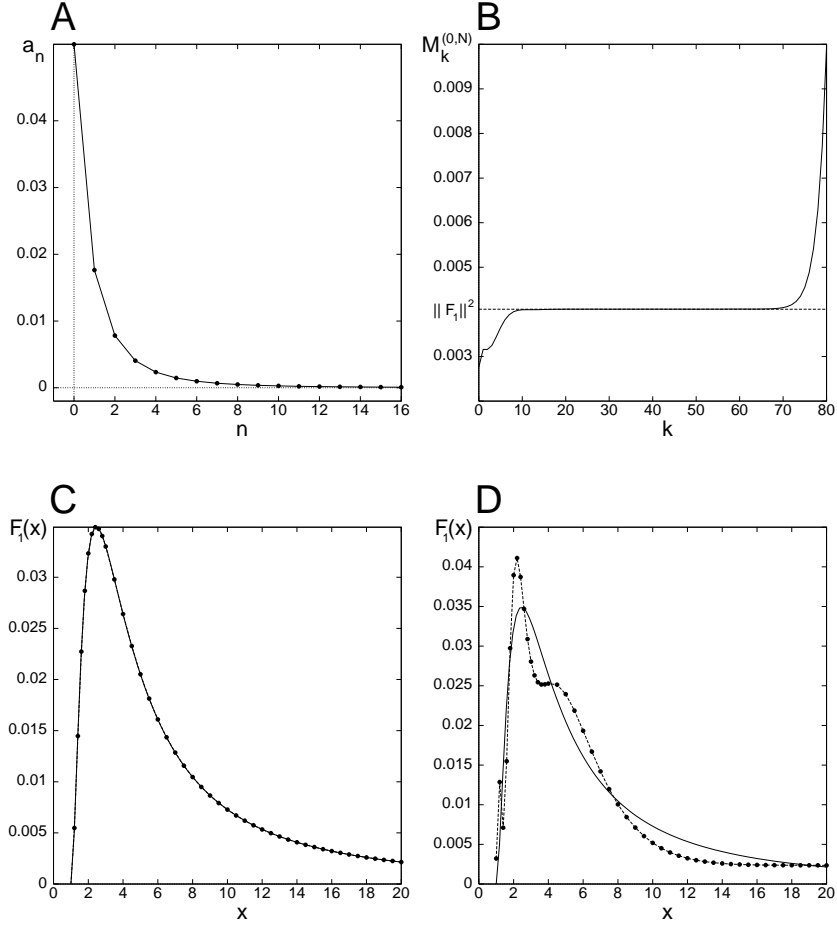}
\caption{\label{figura_1}
Reconstruction of a smooth jump function. $F_1(x)=(x-1)^2/(x^4+x^3+x^2+x+1)$, $x\in [1,+\infty)$;
the number of coefficients $a_n$ used for the computations is $N=30$.
(A) Noiseless coefficients $a_n$ computed according to Eq. (\protect\ref{num1}).
(B) $M_k^{(0,30)}$ versus $k$. The dashed horizontal line indicates the
squared norm of the true function. The algorithm for recovering the truncation index
gave $k_0(0,30)=61$.
(C) Reconstruction of the jump function. The solid line represents the actual function
$F_1(x)$. The dots are samples of the reconstructed jump function $F_1^{(0,30)}$
computed by using $k_0(0,30)=61$.
(D) Comparison between $F_1(x)$ and $F_1^{(0,30)}$ computed by using $k_0(0,30)=75$.}
\end{center}
\end{figure}

\begin{figure}[ht]
\begin{center}
\includegraphics[width=11cm]{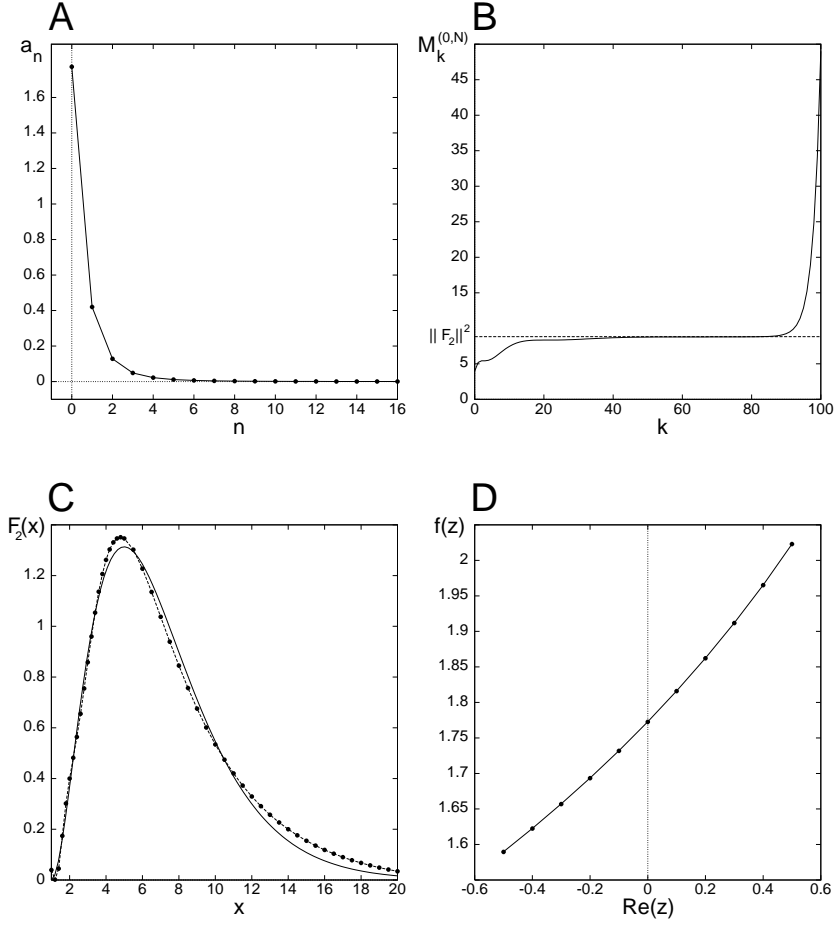}
\caption{\label{figura_2}
Reconstruction of a smooth jump function. $F_2(x)=(x-1)^2\exp(-x/2)$, $x\in [1,+\infty)$; $N=30$.
(A) Noiseless coefficients $a_n$.
(B) $M_k^{(0,30)}$ versus $k$. The dashed horizontal line indicates $\|F_2\|^2$.
The algorithm for recovering the truncation index gave $k_0(0,30)=82$.
(C) Comparison between the actual solution $F_2(x)$ (solid line) and $F_2^{(0,30)}$ (dashed dotted line)
computed by using $k_0(0,30)=82$.
(D) Comparison between the truncated Taylor series $\sum_{n=0}^{30}a_n z^n$ (solid line) and
the approximation computed through the integral of Cauchy type (see Eq. (\protect\ref{sei})) by using the
function $F_2^{(0,30)}$ (dots), for $\Real z \in [-1/2,1/2]$ and $\Imag z = 0$.}
\end{center}
\end{figure}

\begin{figure}[ht]
\begin{center}
\includegraphics[width=11cm]{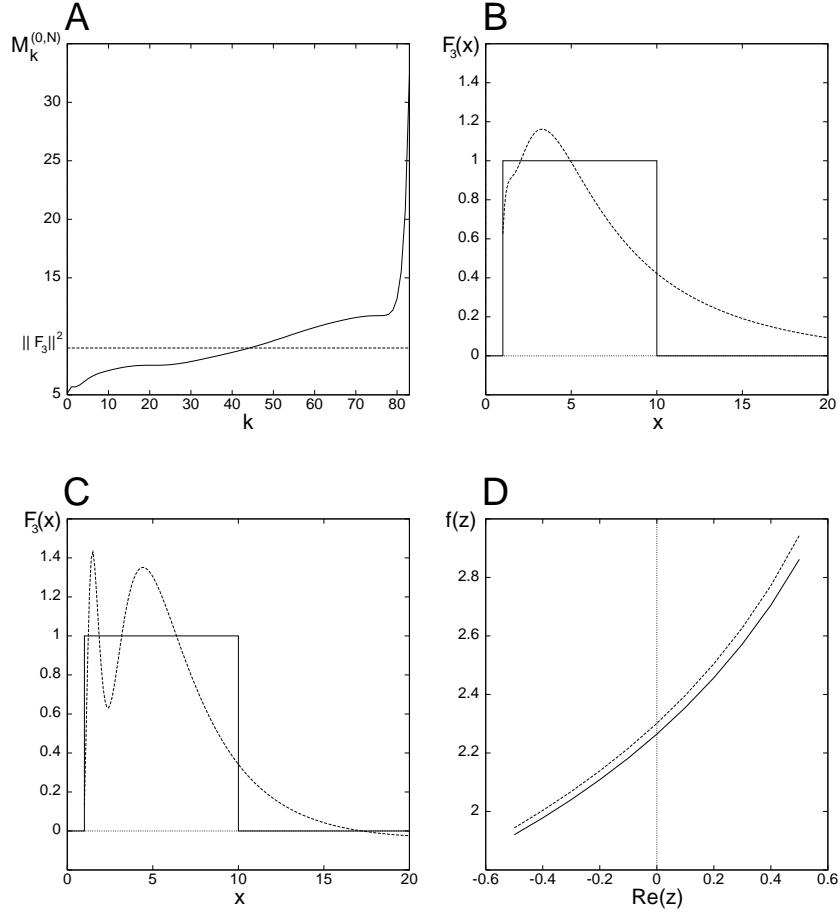}
\caption{\label{figura_3}
Reconstruction of a discontinuous jump function.
$F_3(x)=1$ if $x\in [1,10]$ and $F_3(x)=0$ elsewhere; $N=30$.
(A) $M_k^{(0,30)}$ versus $k$. The dashed horizontal line indicates $\|F_3\|^2$.
The algorithm for recovering the truncation index gave $k_0(0,30)=23$.
(B) Comparison between the actual solution $F_3(x)$ (solid line) and $F_3^{(0,30)}$ (dashed line)
computed by using $k_0(0,30)=23$.
(C) Comparison between the actual solution $F_3(x)$ (solid line) and $F_3^{(0,30)}$ (dashed line)
computed by using $k_0(0,30)=44$.
(D) Comparison between the truncated Taylor series $\sum_{n=0}^{30}a_n z^n$ (solid line) and
the approximation computed through the integral of Cauchy type (see Eq. (\protect\ref{sei})) by using the
function $F_3^{(0,30)}$ with $k_0(0,30)=23$ (dashed line),
for $\Real z \in [-1/2,1/2]$ and $\Imag z = 0$.}
\end{center}
\end{figure}

In this section various numerical aspects of the issues discussed in the
previous section will be illustrated. In particular, the effectiveness
and the accuracy of the reconstruction of the jump function $F(x)$ across the cut will be tested
when only a finite number of coefficients $a_n$ of the Taylor series (\ref{uno}) is known.
Moreover, tests performed with coefficients $a_n$ corrupted by random noise will be
illustrated.
The results hereafter presented summarize a large number of numerical tests performed on
a variety of sample functions, including smooth functions, regular oscillating functions and
discontinuous functions. We shall illustrate only the somehow extreme cases, that is, smooth
(see Figs. \ref{figura_1} and \ref{figura_2}) and discontinuous functions (see Fig. \ref{figura_3}),
since the behavior of the oscillating ones is very similar to the case of the smooth functions.

In Fig. \ref{figura_1} the basic steps for the reconstruction of jump function in the
presence of only round-off noise are shown (see the legend for the numerical details).
In this example the test function is $F_1(x)=\frac{(x-1)^2}{x^4+x^3+x^2+x+1}$ $(x\in [1,+\infty)$).
This function is characterized to be smooth so that we expect the corresponding
reconstruction to be quite satisfactory.
Fig. \ref{figura_1}A illustrates a sample of the coefficients $a_n$, which are computed as
\beq
\label{num1}
a_n = \int_1^{+\infty} x^{-n-1} F_1(x)\, dx.
\eeq
Then, the coefficients $a_n$ do satisfy condition (\ref{due}) of Theorem \ref{the:1}.
In Fig. \ref{figura_1}B the plot of $M_k^{(0,N)}$ versus $k$
(see Eq. (\ref{duesette})) is shown.
According to the properties summarized in Lemma \ref{le:1} (see also the remark after the proof of Lemma \ref{le:1}),
a quite extended plateau is clearly present before $M_k^{(0,N)}$ starts increasing rapidly (see Eq. (\ref{dueotto})).
This fact allows us to determine the truncation point $k_0(0,N)$ which is needed
for constructing the regularized approximation $F_1^{(0,N)}$ (see Eq. (\ref{dueventisei})).
A simple algorithm for the automatic determination of the index $k_0$ has been implemented;
it exploits both the properties provided by Lemma \ref{le:1} that $k_0$ lies approximately on a plateau,
and, moreover, that such a plateau must be located before $M_k^{(0,N)}$ starts growing as a power
of $2N$ (see formula (\ref{dueotto})). \\
In the general case with $\epsilon \neq 0$, the knowledge of the asymptotic behavior of
$M_k^{(\epsilon,N)}$ allows for restricting the range of $k_0$ by defining an upper limit $k_a$
$(k_0 < k_a)$, that represents the value of $k$ where approximately the asymptotic behavior
sets in; in practice, $k_a$ is set as the value of $k$ where $M_k^{(\epsilon,N)}$ and its
asymptotic behavior starts being close enough.
Then, the candidate plateaux are located by selecting the extended intervals of $k < k_a$
where the modulus of the first numerical derivative of $M_k^{(\epsilon,N)}$ is sufficiently small.
Finally, $k_0$ is chosen as the largest value of $k$ belonging to the interval which is closest, but
inferior, to $k_a$.

Figure \ref{figura_1}B shows a quite simple situation, where only one plateau is
present in the interval from about $k=15$ to $k=65$, and whose value, as expected from Eq. (\ref{duecinque}),
is approximately the squared norm of the function $F_1(x)$ (see the horizontal dashed line).
Panels C and D show two reconstructions of the function $F_1(x)$; in Fig. \ref{figura_1}C we have $k_0$ = 61,
which is well within the plateau of Fig. \ref{figura_1}B, and the comparison between the actual
function $F_1(x)$ and its approximation $F^{(0,30)}(x)$ (dotted line) shows the good quality of the
reconstruction. Similar results (not displayed), obtained for different values of $k_0$ ($15 < k_0 < 65$),
indicate that the actual choice of $k_0$ is not critical, provided that it is located in the interval corresponding
to the plateau. Figure \ref{figura_1}D shows that
at $k_0 = 75$, that is just after the correct plateau, the reconstruction deteriorates evidently.

Figure \ref{figura_2} illustrates another example of a smooth sample function that gives rise to a
good reconstruction; in this case $F_2(x)=(x-1)^2\exp(-x/2)$ $(x\in [1,+\infty))$. In Fig. \ref{figura_2}D
we compare the approximation of $f(z)$ given by $f_N(z)=\sum_{n=0}^N a_n z^n$ $(\Real z \in [-1/2,1/2],\,\Imag z = 0)$
with the approximation of $f(z)$ obtained by the Cauchy integral (see Eq. (\ref{sei})) representing
the jump function $F_2(x)$ with its regularized approximation $F_2^{(0,30)}(x)$.

For the sake of completeness, it should be mentioned that
the erratic behavior of the noise can rarely produce very short plateaux located
between the true value of $k_0$ and before $M_k^{(\epsilon,N)}$
starts following its asymptotic behavior (i.e. for $k \simeq k_a$).
In this case our procedure could fail to recover the correct value of $k_0$,
and this minor drawback has been solved heuristically by simply rejecting
the plateaux shorter than a given threshold $L=5$.

However, the situation can become more difficult, and an example is given in Fig. \ref{figura_3}.
Here the jump function $F_3(x)$ is a discontinuous one, being constant in the range
$[1,10]$ and null elsewhere (see Fig. \ref{figura_3}B). In this case the lack of smoothness of $F_3$ causes
the plateau of $M_k^{(\epsilon,N)}$ to be hardly visible, and, consequently,
the definition of the truncation point $k_0$ becomes uncertain (see Fig. \ref{figura_3}A).
In this case, our algorithm gave $k_0=23$, which belongs to the only detectable plateau,
and the corresponding reconstruction is shown in Fig. \ref{figura_3}B.
Fig. \ref{figura_3}C displays the reconstruction of the jump function obtained with $k_0=44$,
that represents the index $k$ such that $M_k^{(\epsilon,N)} \simeq \|F_3\|^2$
(see formula (\ref{duecinque})).
In spite of the fact that, according to definition (\ref{duecinque}), $k_0=44$ is the true truncation index,
the corresponding approximation is evidently worse than the one obtained with $k_0=23$ (see Fig. \ref{figura_3}B).
This fact can be easily explained by recalling that our reconstruction strategy rests on a finite
truncation of a series expansion ($k_0 < k_a < \infty$ for $N<\infty$),
and that the convergence to the true function is only asymptotic. Then,
we pay entirely the penalty of the Gibbs-like phenomenon, so that
our procedure can lead to inaccurate reconstructions in the case of discontinuous jump functions.

\begin{figure}[ht]
\begin{center}
\includegraphics[width=11cm]{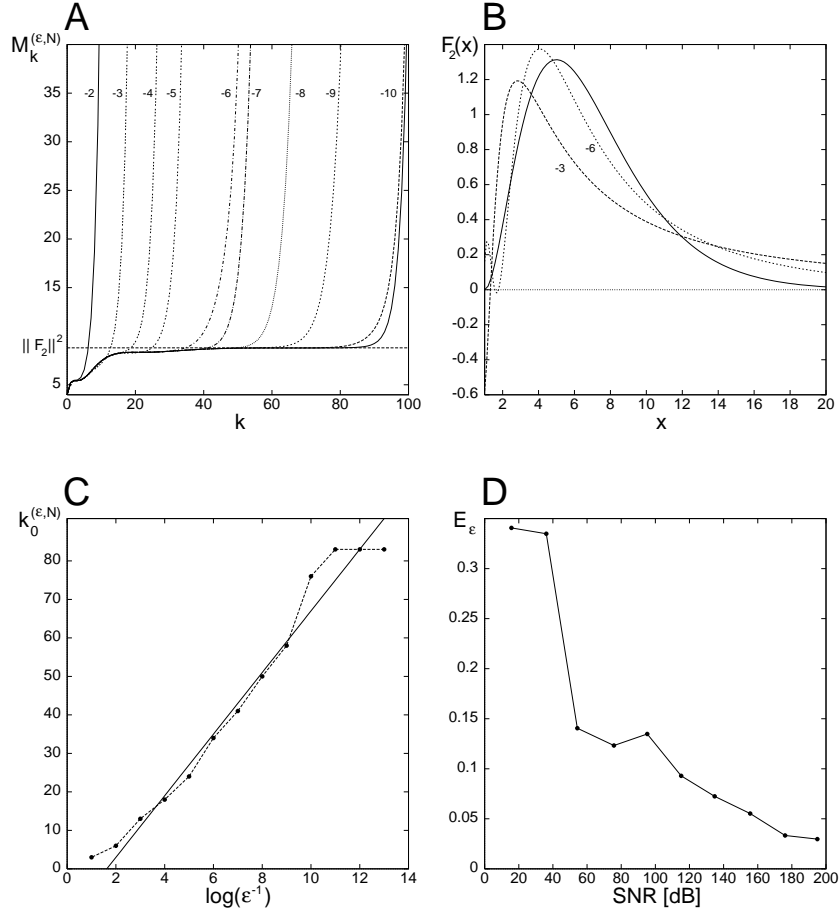}
\caption{\label{figura_4}
Reconstruction of a smooth jump function from noisy coefficients $a_n^{(\epsilon)}$.
The sample function is $F_2(x)$; $N=30$.
The coefficients $a_n$ have been noised by adding white noise uniformly distributed in
the interval $[-\epsilon,\epsilon]$.
(A) $M_k^{(\epsilon,30)}$ versus $k$ computed for $\epsilon$ ranging from
$\epsilon=10^{-2}$ through $\epsilon=10^{-10}$ with step $10^{-1}$. The rightmost solid line
indicates the noiseless $M_k^{(0,30)}$.
(B) Comparison between the actual solution $F_2(x)$ (solid line) and the reconstructions
$F_2^{(10^{-6},30)}$ and $F_2^{(10^{-3},30)}$ (dashed lines).
(C) Plot of $k_0(\epsilon,30)$, set by hand according to its definition (\ref{duecinque}),
against $\log(\epsilon^{-1})$. Each $k_0(\epsilon,30)$ was obtained
as the average over 300 realizations of the corresponding $M_k^{(\epsilon,30)}$.
(D) Plot of the reconstruction error $E_\epsilon$ versus the signal--to--noise ratio SNR.
$E_\epsilon$ is defined as the mean square error of the
reconstruction $F_2^{(\epsilon,N)}(x)$ with respect to the actual function $F_2(x)$, computed
in the interval $x \in [1,20]$. The signal--to--noise ratio is defined
as the ratio of the mean power of the sequence of noiseless coefficients $\{a_n\}$
to the noise variance.}
\end{center}
\end{figure}

An example of the dependence of the reconstruction on the noise is sketched in Fig. \ref{figura_4}. Here the test
function is again $F_2(x)$, and the data illustrated in this figure are representative of an
extensive numerical analysis performed over several different test functions that produced similar results.
The coefficients $a_n$ have been noised by adding white noise,
simulated by computer generated random numbers uniformly distributed in the interval
$[-\epsilon,\epsilon]$. Fig. \ref{figura_4}A groups the functions $M_k^{(\epsilon,N)}$ computed
with different values of $\epsilon$ in the range from $\epsilon=10^{-10}$ through
$\epsilon=10^{-2}$. It is clear that the plateaux become shorter, and, correspondingly,
$k_0$ decreases as $\epsilon$ increases. This fact indicates how, when the noise
increases, the number of coefficients $c_m^{(\epsilon,N)}$ that carry reliable information
for reconstructing the jump function decreases. Examples of such a reconstructions are
plotted in Fig. \ref{figura_4}B. A crude measure of the dependence of the reconstruction stability
on the noise $\epsilon$ is shown in Fig. \ref{figura_4}C, where the truncation index $k_0(\epsilon,N)$,
chosen in this case by hand according to its definition (\ref{duecinque}), is plotted versus the
noise bound $\epsilon$. The plot indicates that the stability estimate is
only logarithmic, as expected for the
problem of the analytic continuation up to the boundary (see \cite{John,Miller}).
However, in spite of this extremely poor stability,
nevertheless it is still possible to obtain acceptable reconstructions of the jump function,
as shown in Fig. \ref{figura_4}B.
Of course, this is not always the case (see for instance the example illustrated in Fig. \ref{figura_3}),
since the quality of the reconstruction depends strongly on the regularity of the jump function,
but, at least for some classes of functions, approximation (\ref{dueventisei}) yields useful results even
in the presence of quite noisy data.
This numerical examples make milder some radical statements that can be found in the literature (see
\cite{John,Miller}) discouraging any attempt of numerical analytic continuation up to the boundary
in physical situations, and indicate that the numerical procedure here proposed
can be successfully applied whenever the function being reconstructed is sufficiently regular.
To this purpose it is worth recalling statement (3) of Theorem \ref{the:1}.
In Fig. \ref{figura_4}D the plot of the reconstruction error, defined as the mean square error
of $F^{(\epsilon,N)}(x)$ with respect to $F(x)$ against the global signal--to--noise ratio (SNR)
is shown. It can be seen that the reconstruction error is still acceptable for SNR up to about 50 dB, whereas
it becomes quite large for smaller SNR.

\section*{Appendix}
\label{appendix_section}
In this Appendix we prove the following lemma on the asymptotic behavior of the
Pollaczek polynomials (see also \cite{Ismail}).

\begin{lemma}
\label{le:appendix}
The Pollaczek polynomials $P_m[-i(n+1/2)]$ satisfy the following
asymptotic behavior for large values of $m$ (at fixed $n$),
\beq
\label{auno}
P_m\left [-i\left (n+\frac{1}{2}\right )\right ] \sim \frac{(-1)^m i^m}{n!} (2m)^n.
\eeq
\end{lemma}

\begin{proof}
Let us introduce the following functions $Q_m^{(1)}(y)$ and $Q_m^{(2)}(y)$, that
we call the associate Pollaczek functions:
\beq
\label{adue}
Q_m^{(1)}(y) = \frac{i}{2}\int_{-\infty}^{+\infty}
\frac{P_m(x) w(x)}{x-y}\,dx \qquad (\Imag y < 0),
\eeq
and
\beq
\label{atre}
Q_m^{(2)}(y) = \frac{i}{2}\int_{-\infty}^{+\infty}
\frac{P_m(x) w(x)}{x-y}\,dx \qquad (\Imag y > 0),
\eeq
where $w(x)=\frac{1}{\pi} \Gamma(1/2+ix)\Gamma(1/2-ix)$. Let us now consider the
following function $H(t,y)$ defined by
\beq
\label{aquattro}
H(t,y) = \frac{i}{2}\int_{-\infty}^{+\infty}
\frac{\left\{\sum_{m=0}^\infty t^m P_m(x)\right\} w(x)}{x-y}\,dx \qquad (\Imag y < 0;\,|t|<1).
\eeq
The series $\sum_{m=0}^\infty t^m P_m(x)$ converges to the generating function
of the Pollaczek polynomials for $|t|<1$ \cite{Bateman2,Szego}. The generating function is given by:
$G(t,x)=(1-it)^{ix-1/2}(1+it)^{-ix-1/2}$. Then, we can rewrite integral (\ref{aquattro})
as
\begin{eqnarray}
\label{acinque}
H(t,y) = \frac{i}{2\pi}\int_{-\infty}^{+\infty}
\frac{(1-it)^{ix-1/2}(1+it)^{-ix-1/2}\left\{\Gamma(\frac{1}{2}+ix)\Gamma(\frac{1}{2}-ix)\right\}}
{x-y}\,dx \nonumber \\
(\Imag y < 0;\,|t|<1).
\end{eqnarray}
We can then integrate the r.h.s. of (\ref{acinque}) by the contour integration method.
Let us, indeed, consider the following integral for $\Real t=0, 0<\Imag t< 1$ and $\Imag y<0$,
\beq
\label{asei}
I(t,y)=\frac{1}{2\pi} \oint_{\gamma_1}
\left [\frac{(1-it)}{(1+it)}\right ]^z
\frac{[(1-it)(1+it)]^{-1/2}}{(-iz-y)}
\left\{\Gamma\left (\frac{1}{2}+z\right )\Gamma\left (\frac{1}{2}-z\right )\right\}\,dz,
\eeq
where $\gamma_1$ is the path shown in Fig. \ref{figura_5}.

\begin{figure}
\begin{center}
\includegraphics[width=7cm]{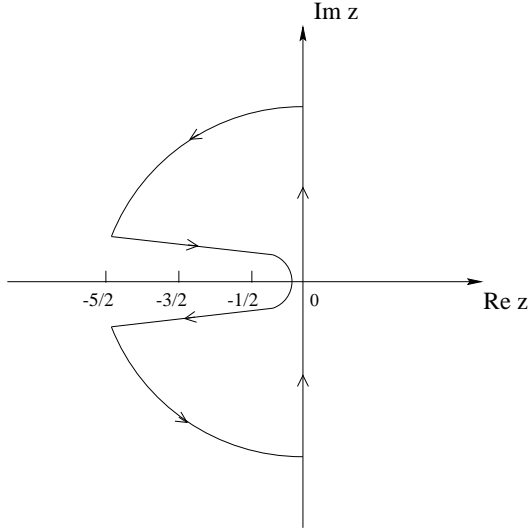}
\caption{\label{figura_5}
Contour used for evaluating integral (\protect\ref{asei}).}
\end{center}
\end{figure}

Rewriting the
term $[(1-it)/(1+it)]^z$ as $\exp (z\ln\frac{(1-it)}{(1+it)})$, it can be easily checked that the contributions
along the quarters of circle belonging to the path $\gamma_1$ vanishes
as $\Real z \rightarrow -\infty$. Therefore, through the Cauchy theorem applied to integral (\ref{asei}),
we obtain, for $\Real t = 0, 0<\Imag t<1, \Imag y < 0$:
\beq
\label{asette}
H(t,y)=\frac{1}{(1-it)}\sum_{k=0}^\infty
\frac{1}{(k+iy+1/2)}\left (\frac{it+1}{it-1}\right )^k,
\eeq
in view of the fact that the only singularities which play a role are the simple poles of
$\Gamma(z+1/2)$ at $z=-k-1/2$. Let us now recall that
\beq
\label{aotto}
_2F_1(1,a,a+1;z)=\sum_{k=0}^\infty\frac{a}{(a+k)}z^k \qquad (|z|<1,\,a\in\C),
\eeq
where $_2F_1$ denotes the Gauss hypergeometric function; we can thus rewrite the expansion (\ref{asette})
as
\beq
\label{anove}
H(t,y)=\frac{1}{(1-it)(iy+1/2)}\,\, \null_2F_1\left (1,iy+\frac{1}{2}, iy+\frac{3}{2};
\frac{it+1}{it-1}\right ),
\eeq
which shows that $H(t,y)$ is analytic in the half--plane $\Imag t > 0$.

By the use of the following relationship, which holds true for the Gauss hypergeometric
function \cite{Bateman1}:
\beq
\label{adieci}
_2F_1(a,b,c;z)=(1-z)^{-a} \,_2F_1 \left (a,c-b,c;\frac{z}{z-1}\right )\qquad
\left (|z|<1, \left |\frac{z}{z-1}\right | < 1\right ),
\eeq
we can rewrite the r.h.s. of formula (\ref{anove}) as
\beq
\label{aundici}
H(t,y)=\frac{1}{2(iy+\frac{1}{2})}~_2F_1\left (1,1,iy+\frac{3}{2}; \frac{1+it}{2}\right ).
\eeq
The hypergeometric series at the r.h.s. of formula (\ref{aundici}) converges inside
a circle of radius 2 and center $t=i$, belonging to the complex $t$-plane. We can then
calculate the following derivatives:
$\frac{1}{m!}(d^mH(t,y)/dt^m)_{t=0}$. By the use of the standard formula for the derivative
of the hypergeometric function we obtain \cite{Bateman1}
\beq
\label{adodici}
\begin{split}
\frac{1}{m!}\left (\frac{d^m H(t,y)}{dt^m}\right )_{t=0} &=
\frac{i^m}{2^{m+1}}\frac{\Gamma(m+1)\Gamma(iy+\frac{1}{2})}{\Gamma(m+iy+\frac{3}{2})} \\
& \times ~_2F_1 \left (m+1, m+1,m+iy+\frac{3}{2};\frac{1}{2}\right ).
\end{split}
\eeq
Next, returning to the integral representation of $H(t,y)$ given by formula (\ref{aquattro}),
we exchange the integral with the sum. This is legitimate for $\Real t=0$, $|t|<1$, since
$|(\frac{1-it}{1+it})^{ix}|=|\exp[ix\ln\frac{1-it}{1+it}]| = 1$, $(x\in\R)$; then we write,
for $\Real t = 0$, $|t|<1$ and $\Imag y < 0$:
\beq
\label{atredici}
H(t,y)=\sum_{m=0}^\infty t^m\left\{\frac{i}{2}\int_{-\infty}^{+\infty}
\frac{P_m(x)w(x)}{x-y}\,dx\right\}=\sum_{m=0}^\infty t^m Q_m^{(1)}(y),
\eeq
and therefore
\beq
\label{aquattordici}
\begin{split}
Q_m^{(1)}(y) &= \frac{1}{m!}\left (\frac{d^m H^1(t,y)}{dt^m}\right )_{t=0} =
\frac{i^m}{2^{m+1}}\frac{\Gamma(m+1)\Gamma(iy+\frac{1}{2})}{\Gamma(m+iy+\frac{3}{2})} \\
& \times ~_2F_1 \left (m+1, m+1,m+iy+\frac{3}{2};\frac{1}{2}\right ).
\end{split}
\eeq
Let us now observe that the function $~_2F_1(a,b,c;z)/\Gamma(c)$
is an entire function of the variable $c$; therefore $Q_m^{(1)}(y)$ can be analytically continued
in the complex $y$-plane, through formula (\ref{aquattordici}) and the only singularities
that it presents in this (open) complex plane are the simple poles of the gamma function
$\Gamma(iy+1/2)$: i.e. $y=i(N+1/2)$, $(N=0,1,2,\ldots)$. Now we can employ the following
relationship which holds true for the hypergeometric function:
\beq
\label{aquindici}
_2F_1(a,b,c;z)=(1-z)^{c-a-b}~_2F_1(c-a,c-b,c;z) \qquad (|z|<1),
\eeq
and we arrive at
\beq
\label{asedici}
Q_m^{(1)}(y)=i^m 2^{-iy-1/2} \frac{\Gamma(m+1)\Gamma(iy+\frac{1}{2})}{\Gamma(m+iy+\frac{3}{2})}
~~_2F_1 \left (iy+\frac{1}{2}, iy+\frac{1}{2},m+iy+\frac{3}{2};\frac{1}{2} \right ).
\eeq

Finally, recalling the asymptotic behavior of
$_2F_1(a,b,c;z)$ for large $|c|$, at fixed $a$, $b$, $z$ (see \cite[formula 10, p. 76]{Bateman1}),
and the asymptotic behavior of the gamma function, we obtain
\beq
\label{adiciassette}
Q_m^{(1)}(y) \sim i^m \Gamma\left (iy+\frac{1}{2}\right ) (2m)^{-iy-1/2}.
\eeq

Let us now focus our attention on the function $Q_m^{(2)}(y)$, defined for $\Imag y > 0$.
If we change $x$ into $-x$, by observing that
\begin{subequations}
\label{adiciotto}
\begin{eqnarray}
&& P_m(-x) = (-1)^m P_m(x), \label{adiciotto.a} \\
&& w(-x)=w(x), \label{adiciotto.b}
\end{eqnarray}
\end{subequations}
we obtain
\beq
\label{adiciannove}
Q_m^{(2)}(y)=(-1)^{m+1}\,\frac{i}{2}\int_{-\infty}^{+\infty}
\frac{P_m(x)w(x)}{x+y}\,dx = (-1)^{m+1} Q_m^{(1)}(-y).
\eeq
This latter relationship, proved for $\Imag y >0$, can be extended by analytic continuation.
Next, by the use of the Cauchy formula we obtain
\beq
\label{aventi}
Q_m^{(2)}(y)-Q_m^{(1)}(y)=\frac{i}{2}\oint_{\gamma_2}\frac{P_m(x)w(x)}{x-y}\,dx=-\pi w(y) P_m(y),
\eeq
where $\gamma_2$ is a path which encircles the singularity located at $x=y$ in a
counterclockwise sense.
From (\ref{adiciannove}) and (\ref{aventi})
we finally obtain
\beq
\label{aventuno}
\pi w(y)P_m(y)=Q_m^{(1)}(y)+(-1)^m Q_m^{(1)}(-y).
\eeq
From the asymptotic behavior (\ref{adiciassette}) and by the use of (\ref{aventuno})
we obtain:
\beq
\label{aventidue}
\begin{split}
P_m(y) &\sim \frac{i^m}{\Gamma(\frac{1}{2}+iy)\Gamma(\frac{1}{2}-iy)}
\left\{\Gamma\left(\frac{1}{2}+iy\right)(2m)^{-iy-1/2}\right . \\
&+ \left . (-1)^m \Gamma\left(\frac{1}{2}-iy\right)(2m)^{iy-1/2} \right\}.
\end{split}
\eeq
Finally, by putting $y=-i(n+1/2)$, we get:
\beq
\label{aventitre}
\begin{split}
P_m\left[-i\left(n+\frac{1}{2}\right)\right]&\sim
i^m \left\{\frac{1}{\Gamma(-n)}(2m)^{-n}+(-1)^m\frac{1}{\Gamma(n+1)}(2m)^n\right\} \\
&= \frac{i^m(-1)^m}{n!}(2m)^n,
\end{split}
\eeq
which is exactly the result we want to prove.
\end{proof}

\end{document}